\documentclass[12pt,twoside]{amsart}

\newtheorem{Thm}{Theorem}

\newtheorem{Lemma}[Thm]{Lemma}
\newtheorem{Prop}[Thm]{Proposition}
\newtheorem*{Prop*}{Proposition}

\theoremstyle{definition}

\newtheorem{Ex}{Example}

\input{diagrams}

\title{K-homology of the rotation algebras $A_{\theta}$}
\author[T.~Hadfield]{Tom Hadfield}
\address{Department of Mathematics, University of California, Berkeley, CA 94720}
\email{hadfield@math.berkeley.edu}
\date{\today}

\begin{document}

\begin{abstract} We study the K-homology of the rotation algebras $A_{\theta}$ using the six term cyclic sequence for the K-homology of a crossed product by ${\bf Z}$. In the case where $\theta$ is irrational, we use Pimsner and Voiculescu's work on AF-embeddings of the $A_{\theta}$ to search for the missing generator of the even K-homology. 
\end{abstract}

\maketitle

\section{Introduction}

In this paper we are concerned with the two dimensional noncommutative tori, the rotation algebras $A_{\theta}$. For $\theta \in$ $[0,1)$, we define $A_{\theta}$ to be the universal C*-algebra generated by unitaries $U$ and $V$ satisfying the relation $VU = \lambda UV$, where $\lambda =$ $e^{2 \pi i \theta}$. 
 These algebras have been extensively studied from many different viewpoints. 
 A thorough overview of the literature appears in Rieffel's survey article \cite{rieffel90}. 

We study the K-homology of the rotation algebras, by which we mean the Kasparov groups ${KK^i}(A_{\theta},{\bf C})$ ($i = 0,1$), and in particular we are interested in exhibiting the generating Fredholm modules. 
 We make extensive use of the six term cyclic sequence for K-homology of a crossed product by ${\bf Z}$, dual to the Pimsner-Voiculescu sequence \cite{pv80} on K-theory. 
 In the commutative situation $\theta =0$, all four generators of the 
K-homology can be exhibited concretely. Three of these Fredholm modules generalize immediately to the case where $\theta \neq 0$,  
 however the canonical ``zero dimensional'' Fredholm module ${\bf z}_0$ 
vanishes. 
In  the final section of this paper we attempt to describe this missing generator, via Pimsner and Voiculescu's work on embedding the $A_{\theta}$ in AF-algebras \cite{pv80a}. 

We note that the K-homology of the rotation algebras was previously studied by Popa and Rieffel in \cite{pr}, who calculated the Ext groups. 
 This approach predated the formalism of Fredholm modules.

\section{Fredholm modules as K-homology}

 Recall that a Fredholm module over a *-algebra 
 $A$ 
is a triple 
 $(H,\pi,F)$,
 where $\pi$ is a\\
  *-representation of $A$ as bounded operators on the Hilbert space $H$. 
 The operator 
$F$ is a selfadjoint element of  
${\bf B}(H)$,
satisfying $F^2 =1$, 
such that the commutators 
$[F,\pi(a)]$ are compact operators for all 
$a \in A$. 
Such a Fredholm module is called odd.

An even Fredholm module is the above data, 
together with a 
${\bf Z}_2$-grading 
of the Hilbert space $H$,
given by a grading operator 
$\gamma \in {\bf B}(H)$
with 
$\gamma = {\gamma}^{*}$, 
$\gamma^2 =1$,
$[\gamma, \pi(a)]=0$ 
for all 
$a \in A$, 
and
$F\gamma = - \gamma F$. 
In general the *-algebra 
$A$ 
will be a dense subalgebra of a C*-algebra, 
closed under holomorphic functional calculus.
Fredholm modules should be thought of as abstract elliptic operators, 
since they are motivated by axiomatizing the important properties of elliptic pseudodifferential operators on closed manifolds. 

This definition is due to Connes \cite{connes94}, p288.
In Kasparov's framework the K-homology groups are given by specializing the second variable in the KK-functor to be 
the complex numbers ${\bf C}$.
Equivalence classes of
even Fredholm modules make up the even K-homology group ${KK^0}(A,{\bf C})$.
Odd Fredholm modules make up the odd K-homology ${KK^1}(A,{\bf C})$. 
 A Fredholm module is said to be degenerate if 
 $[F, \pi(a)] =0$ for all $a \in A$. 
 Degenerate Fredholm modules represent the identity element of the corresponding K-homology group.

Two simple examples of an even and an odd Fredholm module that we will use extensively in the sequel, are as follows.
 Let $A$ be a C*-algebra, $H$ a finite-dimensional Hilbert space, and $\varphi : A \rightarrow {\bf B}(H)$ a *-homomorphism.

\begin{Ex}
\label{evenfred}
We construct a canonical even Fredholm module
 ${\bf z}_0 \in {KK^0}(A, {\bf C})$ : 
\begin{equation}
{\bf z}_0 = (
H_0 = H \oplus H, \pi_0 = \varphi \oplus 0 , 
F_0 = \left(
\begin{array}{cc}
0 & 1 \cr
1 & 0 \cr
\end{array}
\right),
\gamma = 
\left(
\begin{array}{cc}
1 & 0 \cr
0 & -1 \cr
\end{array}
\right)
).
\end{equation}
 In general ${\bf z}_0$ may well represent a trivial element of the even K-homology of $A$ (for example, if $\varphi$ is the zero homomorphism.) 
However, if $A$ is unital, and $\varphi$ is a nonzero  *-homomorphism, the Chern character of ${\bf z}_0$ pairs nontrivially with $[1] \in {K_0}(A)$, 
 showing that ${\bf z}_0$ is a  nontrivial element of K-homology, 
and also that $[1] \neq 0 \in {K_0}(A)$. More precisely :

\begin{Lemma} 
\label{canonicalevenpairing}
Suppose $e = (e_{ij}) \in {M_q}(A)$ is a projection, 
$e = e^{*} = e^2$. 
 Then  
 $< {ch}_{*}( {\bf z}_0 ), [e]> = {\Sigma_{k=1}^q} Tr( \varphi(e_{kk}))$.
 In particular, if $A$ is unital, and $\varphi$ nonzero, then 
 $< {ch}_{*}( {\bf z}_0 ), [1]>$ $=1.$
\end{Lemma}

\begin{proof} 
 Here, 
${ch}_{*} : {KK^0}(A,{\bf C}) \rightarrow {HC^{even}}(A)$, 
is the even Chern character as defined in \cite{connes94}, p295, mapping the even K-homology of $A$ into even periodic cyclic cohomology,
and $<.,.>$ denotes the pairing between K-theory and periodic cyclic cohomology defined in \cite{connes94}, p224.
 We give all the details of this calculation, to avoid later repetition.
We have
\begin{equation}
 < {ch_{*}}( {\bf z}_0 ), [e] > = 
{{\lim}_{n \rightarrow \infty} } {(n!)}^{-1} 
{\Sigma_{i_0, i_1, ..., i_{2n} =1}^{q}}
\psi_{2n}  (e_{i_0 ,i_1}, e_{i_1 ,i_2},...,e_{i_{2n}, i_0}),
\end{equation}
 where (for each $n$) $\psi_{2n}$ is the cyclic $2n$-cocycle defined by
\begin{equation}
 \psi_{2n} ( a_0, a_1, ..., a_{2n}) = 
(-1)^{n(2n-1)} \Gamma (n+1) 
 Tr( \gamma \pi_0 (a_0) [ F_0 , \pi_0 (a_1)] ... [F_0 , \pi_0 (a_{2n})]). 
\end{equation}
 Since $\Gamma(n+1) = n!$ it follows that
\begin{equation}
< {ch_{*}}( {\bf z}_0 ), [1] > = 
{{\lim}_{n \rightarrow \infty} }
 {(-1)}^n
{\Sigma_{i_0, i_1, ..., i_{2n} =1}^{q}}
 Tr( \gamma \pi_0 (e_{i_0 ,i_1}) [F_0 , \pi_0 ( e_{i_1, i_2})] ...
[F_0 , \pi_0 ( e_{i_{2n}, i_0})] ).
\end{equation} 
 Now, for any $a \in A$, 
$[F_0 , \pi_0 (a) ] = 
\left(
\begin{array}{cc}
0 & - \varphi(a) \cr
\varphi(a) & 0 \cr
\end{array}
\right),$
 hence
\begin{equation*}
{\Sigma_{i_0, i_1, ..., i_{2n} =1}^{q}}
 Tr( \gamma \pi_0 (e_{i_0 ,i_1}) [F_0 , \pi_0 ( e_{i_1, i_2})] ...
[F_0 , \pi_0 ( e_{i_{2n}, i_0})] )
\end{equation*}
\begin{equation}
 = {\Sigma_{i_0 =1}^q}
(-1)^n 
 Tr (
\left(
\begin{array}{cc}
\varphi( e_{i_0 ,i_0})  & 0 \cr
0 & 0 \cr
\end{array}
\right)) =
{(-1)}^n {\Sigma_{k=1}^q} Tr (\varphi ( e_{k,k})). 
 \end{equation}
 Therefore
\begin{equation} 
< {ch_{*}}( {\bf z}_0 ), [e] > = 
 {\Sigma_{k=1}^{q}} Tr (\varphi ( e_{k,k}))
\end{equation}
as claimed. If $A$ is unital, and $\varphi$ nonzero, then $\varphi(1) =1$, which proves the second claim. 
\end{proof}

\end{Ex}

\begin{Ex}
\label{oddfred}
We also describe a canonical odd Fredholm module
 ${\bf z}_1 \in {KK^1}( A \times_{\alpha} {\bf Z}, {\bf C})$ : 
\begin{equation}
 {\bf z}_1 = (H_1 = {l^2}( {\bf Z}, H), \pi_1 , F_1)
\end{equation}
 Take 
 $\pi_1 : A \times_{\alpha} {\bf Z} \rightarrow {\bf B}( {l^2}({\bf Z}, H))$ 
to be defined by 
\begin{equation}
(\pi_1 (a) \xi)(n)=\varphi( {\alpha^{-n}}(a)) \xi(n),
\quad
 (\pi_1 (V) \xi)(n)= \xi(n-1),
\end{equation}  
 for $\xi \in {l^2}({\bf Z}, H)$, $a \in A$, and $V$ the unitary implementing the action of ${\bf Z}$ on $A$ 
(via $Va{V^{*}} = \alpha (a)$).  
 Then $\pi_1$ is the usual representation of $A \times_{\alpha} {\bf Z}$ induced from the representation $\varphi$ of $A$.
 We take 
\begin{equation}
F_1 \xi(n) = sign(n) \xi(n) =
\left\{
\begin{array}{cc}
 \xi(n) & : n \geq 0 \cr
- \xi(n) & : n<0 \cr
\end{array}
\right. 
\end{equation}
 It is immediate that $[F_1 , \pi_1 (a)] =0$ for all $a \in A$.
 Further,
  $[F_1 , \pi_1(V)]$ is a finite-rank operator and hence compact, 
 provided $H$ is finite-dimensional. 
Nontriviality of ${\bf z}_1$ (even if $\varphi$ is the zero homomorphism) follows from :

\begin{Lemma}
\label{canonicaloddpairing}
$< {ch}_{*} ({\bf z}_1 ), [V]> = dim(H)$. 
\end{Lemma}

\begin{proof} 
 Again, 
${ch}_{*} : {KK^1}(A,{\bf C}) \rightarrow {HC^{odd}}(A)$, 
is the odd Chern character as defined in \cite{connes94}, p296, mapping the odd K-homology of $A$ into odd periodic cyclic cohomology.
 It is straightforward to calculate this pairing directly, however it is quicker to appeal to Connes' index theorem \cite{connes94}, p296, which states that 
\begin{equation}
< {ch_{*}}( {\bf z}_1 ), [V]> 
 = Index (EVE),
\end{equation}
 where $E = {\frac{1}{2}}(1 + F)$ is the natural orthogonal projection
 ${l^2}({\bf Z}, H) \rightarrow {l^2}({\bf N},H)$.
 We have 
\begin{equation}
 Index(EVE) = dim ker(EVE) - dim ker(E{V^{*}}E) = dim H - 0 = dim H,
\end{equation}
 which proves the result. 
 This shows that ${\bf z}_1$ is a nontrivial Fredholm module, and also that $[V] \neq 0 \in$ 
${K_1}(A \times_{\alpha} {\bf Z})$. 
\end{proof}
\end{Ex}

The Fredholm module  ${\bf z}_0$ can be defined more generally, by taking $\varphi : A \rightarrow {\bf K}(H)$, compact operators on a Hilbert space $H$. This will not work for ${\bf z}_1$, since in this situation the commutator $[F,\pi_1 (V)]$ fails to be compact, unless $H$ is finite dimensional.  
 A very useful special case is when we just have $\varphi : A \rightarrow {\bf C}$. 

\section{Six term cyclic sequence for K-homology}

We now consider the six term cyclic exact sequence for K-homology of crossed products by ${\bf Z}$, dual to the Pimsner-Voiculescu sequence for K-theory, as described in \cite{blackadar}, p199.

Recall \cite{pv80} that associated to any crossed product algebra $A \times_{\alpha} {\bf Z}$  
is the following semisplit short exact sequence of C*-algebras, 
the Pimsner-Voiculescu ``Toeplitz extension" 
\begin{equation}
\label{toeplitzext}
 0 \rightarrow {A \otimes {\bf K}} \rightarrow {T_{\alpha}} \rightarrow  {A \times_{\alpha} {\bf Z}} \rightarrow  0.
\end{equation}
 Here ${T_{\alpha}}$ 
is the C*-subalgebra of 
$( {A {\times}_{\alpha} {\bf Z}}) \otimes {\it T} $
generated by 
$a \otimes 1$, 
$a \in A$ 
and 
$ V \otimes f$,
where 
$V$ 
is the unitary implementing the action of 
$\alpha$ 
on 
$A$, 
and 
$f$ 
is the non-unitary isometry generating the ordinary Toeplitz algebra
 $T$, 
that is  
$f \in {\bf B}({l^2}({\bf N}))$,
$f{e_n} = {e_{n+1}}$.
 This extension defines the Toeplitz element 
${\bf x} \in {KK^1}( {A \times_{\alpha} {\bf Z}}, A) $.  

Applying the K-functor gives the Pimsner-Voiculescu six term cyclic sequence for K-theory. 
The corresponding six term cyclic sequence for K-homology is:

\begin{diagram}
{{KK^0}( A  , {\bf C})} & \lTo^{id - {\alpha}^{*}} & {{KK^0}( A, {\bf C})} & \lTo^{i^{*}} & {{KK^0}(A {\times_{\alpha}} {\bf Z}, {\bf C})}  \\
\dTo^{\partial_0} & & & & \uTo^{\partial_1} \\
{{KK^1}( A {\times_{\alpha}} {\bf Z} , {\bf C})} & \rTo^{i^{*}} & {{KK^1}( A, {\bf C})} & \rTo^{id - {\alpha}^{*}} & {{KK^1}( A, {\bf C} )}\\
\end{diagram}

Here $i$ denotes the canonical inclusion map 
 $i : A \hookrightarrow A \times_{\alpha} {\bf Z}$. 
The vertical maps 
${\partial_0}$ and  ${\partial_1}$  
are given by taking the Kasparov product with the  Toeplitz element:
\begin{equation}
\partial_i : {KK^i}(A,{\bf C}) \rightarrow {KK^{i+1}}( A \times_{\alpha} {\bf Z}, {\bf C}), \quad
{\bf z} \mapsto {\bf x} {\hat{\otimes}}_A {\bf z}
\end{equation}
 This sequence formulated in terms of Ext appears in the original paper of Pimsner and Voiculescu \cite{pv80}.
However, the relationship between Ext and the Fredholm module picture of K-homology is not transparent.

Let $A$ be a C*-algebra, with a finite-dimensional representation  $\varphi : A \rightarrow {\bf B}(H)$.  We assume for convenience that $\varphi(A) H = H$. 
Then the  Fredholm modules ${\bf z}_0$ and ${\bf z}_1$ described above (Examples \ref{evenfred} and \ref{oddfred}) are related via the morphism  $\partial_0$ as follows.

\begin{Prop}
\label{bdy0}
Under the map ${\partial_0}$ 
we have $\partial_0 ( {\bf z}_0 ) = {\bf z}_1$. 
\end{Prop}
\begin{proof}
We describe the Pimsner-Voiculescu Toeplitz element
${\bf x} \in {KK^1}( A {\times}_{\alpha} {\bf Z}, A )$
as the Kasparov triple 
\begin{equation}
\label{defntoeplitzelement}
({E_1}, {\phi_1}, {F_1}) \in
 {\bf E}( A {\times}_{\alpha}{\bf Z} , A {\hat{\otimes}} {\bf C}_1 ).
\end{equation}
 Here, ${\bf E}(B,D)$ denotes the set of Kasparov triples over a pair of C*-algebras $B$, $D$, \cite{blackadar}, p143. 
 We take 
${E_1} = {l^2}({\bf Z},A) {\hat{\otimes}} {\bf C}_1$, 
 with the obvious
$A {\hat{\otimes}} {\bf C}_1 $-valued inner product.
 Here $\hat{\otimes}$ is the graded tensor product of Hilbert modules, while ${\bf C}_1$ is the Clifford algebra of the one dimensional complex vector space ${\bf C}$, being generated by elements 1 and $\epsilon$, with $\epsilon^2 =1$. 
We have
$\phi_1 : A {\times}_{\alpha}{\bf Z} \rightarrow {\bf B}({E_1})$
 given by
\begin{equation}
{{\phi_1}(x) ( \xi {\hat{\otimes}} \omega ) ={\phi_1}'(x) \xi {\hat{\otimes}} {\omega}}
\end{equation}
where
\begin{equation}
{({\phi_1}'(a) \xi)(n) = {\alpha^n}(a) \xi(n) }, \quad
{({\phi_1}'(V) \xi)(n) = \xi (n+1)}
\end{equation}
 for 
$a \in A$, 
$\xi \in {l^2}({\bf Z},A)$ 
and
$V$ 
is the unitary implementing 
$\alpha$. 
 The operator $F_1$ is given by 
${F_1} = F \hat{\otimes} \epsilon$,
with 
$F \xi (n) = sign(n) \xi(n)$.
 
The canonical even Fredholm module 
${\bf z}_0$ 
corresponds to the Kasparov triple
\begin{equation}
 (H \oplus H, \phi_0 = \varphi \oplus 0, F_0 =
\left(
\begin{array}{cc}
0 & 1 \cr
1 & 0 \cr
\end{array}
\right)
)
 \in {\bf E}(A, {\bf C}).
\end{equation}
Recall that, given C*-algebras $A$, $B$ and  $D$, and Kasparov triples
 $(E_1, \phi_1, F_1) \in$
 ${\bf E}(A,D)$, and 
$(E_2, \phi_2, F_2) \in$ ${\bf E}(D,B)$, then the product is given \cite{blackadar}, p166 by the triple 
\begin{equation}
\label{defnproduct}
(E = E_1 {\hat{\otimes}}_{\phi_2} E_2, \phi = \phi_1 \hat{\otimes} 1, F) 
 \in 
{\bf E} (A,B)
\end{equation}
 where $F$ is a ``suitable''  combination of $F_1$ and $F_2$. 
 Almost all the difficulties involved in calculating the product lie in finding the correct $F$. 

Our calculation of the product proceeds in three steps. 
 First, the triples representing ${\bf x}$ and ${\bf z}_0$ need to be compatible, in that ${\bf z}_0$ should be represented by an element of 
${\bf E}(A \hat{\otimes} {\bf C}_1,..)$ rather than ${\bf E}(A,..)$. 
 This is achieved via the morphism \cite{blackadar}, p160:
 \begin{equation*}
{\tau_{{\bf C}_1}} : {\bf E}(A,{\bf C}) \rightarrow {\bf E}(A {\hat{\otimes}}{\bf C}_1 , {\bf C}_1 )
\end{equation*}
\begin{equation}
\label{defntaumap}
(E, \phi, F) \mapsto ( E \hat{\otimes} {\bf C}_1, \phi \hat{\otimes} id, F \hat{\otimes} 1).
\end{equation}
 where $(\phi \hat{\otimes} id)(a \hat{\otimes} \omega)( \xi \hat{\otimes} \omega')$ $= \phi(a) \xi \hat{\otimes} \omega \omega'$. 
 
The second step is to calculate the product\\ 
 ${\bf x} {\hat{\otimes}}_{A \hat{\otimes} {\bf C}_1} \tau_{{\bf C}_1} ( {\bf z}_0 )$, following the procedure outlined above (\ref{defnproduct}). 
 This gives us a triple in\\
 ${\bf E}( A \times_{\alpha} {\bf Z}, {\bf C}_1)$.

 The final step is to show that this triple represents 
 the same element of 
 ${KK^1}( A \times_{\alpha} {\bf Z}, {\bf C})$ 
 as ${\bf z}_1$.

{\bf Step One:}
We apply the map
$\tau_{{\bf C}_1}$
to get
\begin{equation}
\tau_{{\bf C}_1} ({\bf z}_0 )= 
( (H \oplus H) {\hat{\otimes}} {\bf C}_1 , {\phi_0} {\hat{\otimes}} id, {F_0} {\hat{\otimes}} 1)= 
({E_2}, {\phi_2}, F_2)
 \in 
{\bf E}( A \hat{\otimes} {\bf C}_1, {\bf C}_1).
\end{equation}

{\bf Step Two:} 
Now we can take the product. 
We have 
\begin{equation}
\partial_0 ({\bf z}_0) = {\bf x} {\hat{\otimes}}_{A \hat{\otimes} {\bf C}_1} {\tau_{{\bf C}_1} }({\bf z}_0 )
 = (E, \phi, F)
\in {\bf E}( A \times_{\alpha}{\bf Z}, {\bf C}_1).
\end{equation}
 Here
$E= {E_1} {\hat{\otimes}}_{ \phi_0 {\hat{\otimes}} 1} {E_2}$, 
 $\phi = \phi_1 \hat{\otimes} 1$, and $F$ is yet to be found.  
As elements of $E \cong$\\
 $({l^2}({\bf Z},A) \hat{\otimes} {\bf C}_1 ) \hat{\otimes} ((H \oplus H) \hat{\otimes} {\bf C}_1 )$, 
we have 
\begin{equation}
{(  {\delta_k}a {\hat{\otimes}} {\omega_1} ) {\hat{\otimes}}
( {\bf v} {\hat{\otimes}} \omega_2 )  \sim  
 ({\delta_k} {\hat{\otimes}} 1 ) {\hat{\otimes}} ( {{\phi_0}(a)}{\bf v} {\hat{\otimes}} {\omega_1} {\omega_2})}.
\end{equation}
 So we can identify $E$ as a submodule of  
${l^2}( {\bf Z}, H \oplus H ) {\hat{\otimes}} {\bf C}_1$ 
 via the morphism 
\begin{equation}
\label{simbded}
(  {\delta_k} a_k {\hat{\otimes}} {\omega_1} ) {\hat{\otimes}}
( {\bf v} {\hat{\otimes}} \omega_2 )
 \mapsto 
\delta_k \phi_0 (a_k) {\bf v} \hat{\otimes} \omega_1 \omega_2.
\end{equation}
 Since by assumption $\varphi(A)H = H$, and 
$\phi_0 (a) =$
$\left(
\begin{array}{cc}
\varphi(a) & 0 \cr
0 & 0 \cr
\end{array}
\right)$, 
 the image of this morphism 
 can be naturally identified with ${l^2}({\bf Z},H) \hat{\otimes} {\bf C}_1$. 
 After this identification,    
$\phi = {\phi_1} {\hat{\otimes}} 1$
acts via 
\begin{equation}
\phi(x) (\xi \hat{\otimes} \omega) = \phi'(x) \xi \hat{\otimes} \omega
\end{equation}
with 
\begin{equation}
(\phi'(a) \xi)(n) =
\varphi ( {\alpha^n}(a) )  \xi(n),\quad  
( \phi'(V) \xi)(n) = \xi(n+1)
\end{equation}
 for 
$a \in A$, 
$\xi \in {l^2}({\bf Z},H)$.

 We use the Connes-Skandalis formalism of connections \cite{blackadar}, p170, to find a suitable $F$. 
By \cite{blackadar} Prop 18.10.1, such an $F$ will be given 
by 
\begin{equation}
\label{suitableF}
F = F_1 \hat{\otimes}  1 + ( (1- {F_1}^2)^{1/2}  \hat{\otimes} 1) G 
\end{equation}
 where $G$ is an $({F_2} \hat{\otimes} 1)$-connection.
 By \cite{blackadar}, Prop 18.3.3, abstractly such a connection $G$ must exist. 
 In this situation, since ${F_1}^2 =1$, we can just take $F = {F_1} \hat{\otimes}  1$, and there is no need to find a concrete $G$. 
 This defines $F$ as an operator on $E$. 
 Under our identification of $E$ with ${l^2}({\bf Z},H) \hat{\otimes} {\bf C}_1$, we have 
\begin{equation}
F ( \xi \hat{\otimes} \omega) = F' \xi \hat{\otimes} \epsilon \omega
\end{equation}
 with $F' \xi (n) = sign(n) \xi(n)$. 
 We note that $E$ was originally defined as a submodule of\\ 
 ${l^2}({\bf Z}, H \oplus H) \hat{\otimes} {\bf C}_1$, 
 and this submodule is invariant under the actions of 
$A \times_{\alpha} {\bf Z}$ and $F$ defined above. 
 Hence identifying $E$ with ${l^2}({\bf Z},H) \hat{\otimes} {\bf C}_1 $ is well-defined. 
 Therefore, we have calculated the product and obtained a triple 
\begin{equation}
\label{product}
 (E, \phi, F) \in {\bf E}( A \times_{\alpha} {\bf Z}, {\bf C}_1)
\end{equation}
 representing $\partial_0 ({\bf z}_0)$.

{\bf Step Three:}
 An odd Fredholm module $(H_1, \pi_1, F_1) \in$
 $KK^1 (A \times_{\alpha} {\bf Z}, {\bf C})$ is represented by the Kasparov triple 
 $(H_1 \hat{\otimes} {\bf C}_1, \pi_1 \hat{\otimes} 1, F_1 \hat{\otimes} \epsilon) \in$ 
 ${\bf E}( A \times_{\alpha} {\bf Z}, {\bf C}_1 )$. 
 It is immediate to see that the Kasparov triple described in (\ref{product}) corresponding to the product 
 $\partial_0 ({\bf z}_0 )$ represents the Fredholm module ${\bf z}_1$ 
exactly as in Example \ref{oddfred}.
 Hence 
 $\partial_0 ({\bf z}_0 ) =$ 
 ${\bf z}_1$ as Fredholm modules.
 This  completes the proof. 
\end{proof}

\section{Application to the rotation algebras $A_{\theta}$}

We now illustrate this work with the example of the rotation algebras $A_{\theta}$.  
 Since the $A_{\theta}$ are deformations of the commutative algebra $A_0 =$ $C({\bf T}^2)$ (the case $\theta = 0$) we consider this case first. 

\begin{Prop} 
We have ${KK^i}(A_0, {\bf C}) \cong {\bf Z}^2$, $i=0,1$. 
\end{Prop}
\begin{proof}
It is well known that the K-groups $K_i ( A_0)$ ($i =0,1$) are both isomorphic to ${\bf Z}^2$. The generators of $K_0 (A_0)$ are $[1]$ and the Bott projector $[Bott]$. The generators of $K_1 (A_0)$ are $[U]$ and $[V]$. 
  Now, it follows from Rosenberg and Schochet's universal coefficient theorem \cite{rs}, \cite{blackadar}, p234, that, for a C*-algebra $A$ whose K-groups are free abelian, then ${KK^i}(A, {\bf C}) \cong$ ${K_i}(A)$ (as abelian groups). 
 Hence the result. 
\end{proof}

We describe the generators of the K-homology. 
 First of all, we have a canonical ``zero dimensional'' even Fredholm module ${\bf z}_0$ (Example \ref{evenfred}) corresponding to the *-homomorphism $\varphi : A_0 \rightarrow {\bf C}$ given by $U, V \mapsto 1$. 
 Via an identical calculation to Lemma \ref{canonicalevenpairing}, we have
\begin{equation}
< {ch_{*}}({\bf z}_0), [1]> = 1 = < {ch_{*}}({\bf z}_0), [Bott]>.
\end{equation} 
 Since the pairings with the generators of K-theory are both 1, it follows from Connes' index theorem \cite{connes94}, p296 that this Fredholm module 
 is a generator of K-homology, in the sense that if ${\bf z}$ is another Fredholm module, with ${\bf z}_0 = n {\bf z}$ for some $n \in {\bf Z}$, then $n = \pm 1$. 

 For the odd K-homology, we first  describe $A_0$ as a crossed product by (a trivial action of) ${\bf Z}$ in two obvious ways, first as ${C^{*}}(V) \times_{id} {\bf Z}$, second as 
${C^{*}}(U) \times_{id} {\bf Z}$, where the trivial action of ${\bf Z}$ is implemented by $U$ and $V$ respectively. 
 We denote the corresponding odd Fredholm modules of Example \ref{oddfred} by 
 ${\bf z}_1$ and ${{\bf z}_1}'$. 
 We have 
${\bf z}_1 = ( {l^2}({\bf Z}),{\pi_1},F)$, 
where
\begin{equation}
\label{azerozone}
{\pi_1}(U){e_k} = {e_{k+1}}, \quad 
{\pi_1}(V) = I,
\end{equation}  
and
${{\bf z}_1}'=( {l^2}({\bf Z}), {\pi_1}' , F)$
with 
\begin{equation}
\label{azerozoneprime}
{\pi_1}'(U) = I, \quad  
{\pi_1}'(V){e_k}={e_{k+1}}
\end{equation}
and in each case 
$F e_k =$ $sign(k) e_k$. 
 These Fredholm modules generate the odd K-homology ${KK^1}(A_0, {\bf C})$. 
 Calculations identical to Lemma \ref{canonicaloddpairing} show that 
\begin{equation}
\label{chernzone}
< {ch_{*}}({\bf z}_1), [U]> = 1, \quad
< {ch_{*}}({\bf z}_1), [V]> = 0,
\end{equation}
\begin{equation}
\label{chernzoneprime}
< {ch_{*}}({{\bf z}_1}'), [U]> = 0, \quad
< {ch_{*}}({{\bf z}_1}'), [V]> = 1.
\end{equation}

The second generator of the even K-homology is the Fredholm module ${\bf Dirac}$, which is the bounded formulation of the Dirac operator on ${\bf T}^2$. 
\begin{equation}
\label{azerodirac}
{\bf Dirac} = (H, \pi, F)
\end{equation} 
where 
$H = {l^2}({\bf Z}^2 ) \oplus {l^2}({\bf Z}^2 )$, with $A_0$ acting on the orthonormal basis ${\{ e_{m,n} \}}_{(m,n) \in {\bf Z}^2}$
 for ${l^2}({\bf Z}^2)$ via 
\begin{equation}
U{e_{m,n}} = e_{m+1,n}, \quad
V{e_{m,n}} = e_{m,n+1}
\end{equation} 
 and we take
\begin{equation}
 \pi(a)=
\left(
\begin{array}{cc}
a & 0 \cr
0 & a \cr
\end{array}
\right), \quad
 F=
\left(
\begin{array}{cc}
0 & {F_0} \cr
{F_0}^{*} & 0 \cr
\end{array}
\right)
\end{equation} 
where 
$F_0$ 
is the diagonal operator defined by 
\begin{equation}
{F_0}{e_{m,n}} = 
\left\{
\begin{array}{cc}
{\frac{m+in}{ {(m^2 + n^2)}^{1/2}}}{e_{m,n}} &: (m,n) \neq (0,0) \cr
{e_{0,0}} & : (m,n) =(0,0) \cr
\end{array}
\right.
\end{equation}

We can use the Baum-Connes assembly map \cite{bch} to identify the K-homology and K-theory of $A_0$.  For a general discrete group $\Gamma$, the assembly map is a homomorphism 
\begin{equation}
\mu : {KK^i}( C_0 (B \Gamma), {\bf C}) \rightarrow {K_i} ( {C_{r}^{*}}( \Gamma ))
\end{equation}
 where $B \Gamma$ is the classifying space of the group $\Gamma$. 
 Now, $A_0 =$ $C({\bf T}^2)$ $={C^{*}}({\bf Z}^2)$, and we have $B {\bf Z}^2 =$ ${\bf T}^2$. In this situation, the assembly map is an isomorphism (since very trivially ${\bf Z}^2$ is amenable) and basically acts as a Fourier transform. 
 We have
\begin{equation}
 \mu : {KK^i} ( A_0 , {\bf C}) \cong K_i (A_0)
\end{equation}
\begin{equation}
{\bf z}_0 \mapsto \pm [1],
\end{equation}
\begin{equation}
{\bf z}_1 \mapsto \pm [U], \quad
{{\bf z}_1}' \mapsto \pm [V],
\end{equation}
\begin{equation}
{\bf Dirac} \mapsto \pm [Bott]
\end{equation}
 Hence for the commutative situation everything is transparent. 

We apply this knowledge to the case $\theta \neq 0$. 

\begin{Prop}
 For $0 \leq \theta \leq 1$, the K-homology groups of the $A_{\theta}$ are
${KK^i}( A_{\theta}, {\bf C}) \cong {\bf Z}^2$.  
\end{Prop}
\begin{proof} This again follows as a corollary of Rosenberg and Schochet's the universal coefficient theorem, since the K-groups of the $A_{\theta}$ are both ${\bf Z}^2$ for all values of $\theta$ \cite{rieffel81}. 
\end{proof}

Three of the four generators of the K-homology of $A_0$ generalize immediately to the case $\theta \neq 0$. 
 The odd K-homology is still generated by Fredholm modules 
${\bf z}_1 = ( {l^2}({\bf Z}),{\pi_1},F)$, 
and
${{\bf z}_1}'=( {l^2}({\bf Z}), {\pi_1}' , F)$
where
\begin{equation}
\label{athetazone}
{\pi_1}(U){e_k} = {e_{k+1}}, \quad 
{\pi_1}(V){e_k} = {\lambda^{k}}{e_k},
\end{equation}  
\begin{equation}
\label{athetazoneprime}
{\pi_1}'(U){e_k}={\lambda^{-k}}{e_k}, \quad  
{\pi_1}'(V){e_k}={e_{k+1}}
\end{equation}
and in each case 
$F e_k =$ $sign(k) e_k$.  
 The pairings of the Chern characters of these Fredholm modules with the generators of $K_1 (A_{\theta})$ are unchanged from (\ref{chernzone}), (\ref{chernzoneprime}). 

The Fredholm module ${\bf Dirac}$ is slightly modified.
\begin{equation}
\label{Athetadirac}
{\bf Dirac} = (H, \pi, F)
\end{equation} 
where 
$H = {L^2}({A_{\theta}},\tau ) \oplus {L^2}( {A_{\theta}}, \tau)$, with 
$\tau$ 
being the canonical trace on 
$A_{\theta}$,
given by  
\begin{equation}
{\tau}({\Sigma {a_{m,n}}{U^m}{V^n}}) = {a_{0,0}}
\end{equation}  
We identify 
${L^2}({A_{\theta}},\tau)$ with
${l^2}({\bf Z}^2 )$
with orthonormal basis 
${\{ e_{m,n} \}}_{(m,n) \in {\bf Z}^2}$,
with 
$A_{\theta}$ acting via
\begin{equation}
U{e_{m,n}} = e_{m+1,n}, \quad
V{e_{m,n}} = {\lambda^{m}}e_{m,n+1}
\end{equation} 
 and we take
\begin{equation}
 \pi(a)=
\left(
\begin{array}{cc}
a & 0 \cr
0 & a \cr
\end{array}
\right), \quad
 F=
\left(
\begin{array}{cc}
0 & {F_0} \cr
{F_0}^{*} & 0 \cr
\end{array}
\right)
\end{equation} 
where 
$F_0$ 
is the diagonal operator defined by 
\begin{equation}
{F_0}{e_{m,n}} = 
\left\{
\begin{array}{cc}
{\frac{m+in}{ {(m^2 + n^2)}^{1/2}}}{e_{m,n}} &: (m,n) \neq (0,0) \cr
{e_{0,0}} & : (m,n) =(0,0) \cr
\end{array}
\right.
\end{equation}

We consider the six term exact sequence on K-homology, in the case where our algebra $A$ is 
$C({\bf T})$, 
thought of as 
${C^{*}}(U)$ 
for some generating unitary 
$U$, 
with 
${\bf Z}$-action 
given by the automorphism 
$\alpha(U) = \lambda U$, 
where 
$\lambda = {e^{2 \pi i \theta}}$.
So 
$A \times_{\alpha} {\bf Z} \cong {A_{\theta}}$, 
the rotation algebra.

\begin{diagram}
{{KK^0}(A, {\bf C})} & \lTo^{id - {\alpha}^{*}} & {{KK^0}( A, {\bf C})} & \lTo^{i^{*}} & {{KK^0}( A_{\theta}, {\bf C})}\\ 
\dTo^{\partial_0} & & & & \uTo^{\partial_1} \\
{{KK^1}( A_{\theta} , {\bf C}) } & \rTo^{i^{*}} & {{KK^1}( A, {\bf C})} & \rTo^{id - {\alpha}^{*}} & {{KK^1}( A, {\bf C})} \\
\end{diagram}

Since 
${K_0}(A)$ 
and 
${K_1}(A)$ 
are both isomorphic to 
${\bf Z}$, 
 generated by $[1]$ and $[U]$ respectively, 
the universal coefficient theorem tells us that the K-homology groups
${KK^0}(A,{\bf C})$ 
and  
${KK^1}(A,{\bf C})$ 
are also both 
${\bf Z}$. 
The generator of 
${KK^0}(A,{\bf C})$ 
is the canonical Fredholm module 
${\bf w}_0$ 
 (Example \ref{evenfred}),
 corresponding to the unital *-homomorphism 
$\varphi : A \rightarrow {\bf C}$
 given by
$U \mapsto 1$.
The generator 
${\bf w}_1$ 
of 
${KK^1}(A,{\bf C})$ 
is the Fredholm module 
$({l^2}({\bf Z}), \pi, F)$,
with 
$U$ 
acting on 
${l^2}({\bf Z})$ 
as the bilateral shift,
$\pi(U) {e_k}={e_{k+1}}$,
and 
$F{e_k}=sign(k){e_k}$.

We saw previously that 
${KK^1}( {A_{\theta}},{\bf C}) \cong {\bf Z}^2 $,
with generators
${\bf z}_1$, ${{\bf z}_1}'$ defined in (\ref{athetazone}), (\ref{athetazoneprime}).  
 We know from Theorem \ref{bdy0} that
$\partial_0 ({\bf w}_0 )={{\bf z}_1}' $.
 The inclusion 
$i : A \hookrightarrow A_{\theta}$
induces maps
\begin{equation}
{ {i^{*}}: {KK^j}( {A_{\theta}},{\bf C}) \rightarrow {KK^j}(A,{\bf C})  \quad
(j=0,1)}.
\end{equation}
 \begin{Lemma}
${i^{*}}({{\bf z}_1}') ={\bf 0}$, ${i^{*}}({\bf z}_1)={\bf w}_1$ and 
$(id - {\alpha}^{*})({\bf w}_1) = {\bf 0}$. 
\end{Lemma} 
\begin{proof} 
We have 
${i^{*}}({{\bf z}_1}') = ({l^2}({\bf Z}) , {\pi_1} \circ i , F)$,
a trivial Fredholm module, since 
$[F, \pi_0 \circ i (U)]=0$. 
 We can also see that  
${i^{*}}({\bf z}_1)={\bf w}_1$.
We have 
${\alpha^{*}}({\bf w}_1)$ is the Fredholm module
$({l^2}({\bf Z}), \pi \circ \alpha, F)$,
with 
 $\pi \circ \alpha (U) {e_k}= \lambda {e_{k+1}}$. 
 Hence the Fredholm modules ${\bf w}_1$ and 
 ${\alpha^{*}}({\bf w}_1)$ are unitarily equivalent, via the unitary 
 $Q e_k = {\lambda}^k e_k$, and therefore represent the same element of K-homology. 
 So
$(id - {\alpha}^{*})({\bf w}_1) = {\bf 0} \in$ ${KK^0}(A,{\bf C})$.
\end{proof}

We can see that these results are in agreement with the Pimsner-Voiculescu sequence for K-theory:

\begin{diagram}
{K_0}(A) & \rTo^{id -{\alpha}_{*}} & {K_0}(A) & \rTo^{i_{*}} & {K_0}( A \times_{\alpha} {\bf Z} ) \\
\uTo^{\delta_1}  & & & & \dTo^{\delta_0} \\
{K_1}( A {\times_{\alpha}} {\bf Z} ) & \lTo^{i_{*}} & {K_1}(A) & \lTo^{id - {\alpha}_{*}} & {K_1}(A) \\
\end{diagram}.

It is well known \cite{rieffel81} that  
${K_1}(A_{\theta}) \cong {\bf Z}^2$, 
generated by 
$[U]$, 
$[V]$, 
and
${K_0}(A_{\theta}) \cong {\bf Z}^2$, 
generated by 
$[1]$, 
$[p]$,
where 
$p$ 
is a Rieffel projection of trace 
$\theta$
 (where $\theta \in (0,1)$). 
We know \cite{ap} that
$\delta_0  [p] = [U]$, 
$\delta_0  [1] =  0$,
and also
$\delta_1 [U]=0$, 
$\delta_1  [V] =[1]$. 
Combining this with 
the mapping of K-homology 
\begin{equation}
 \partial_0  : {KK^0}(A,{\bf C}) \rightarrow {KK^1}( A_{\theta},{\bf C})
\end{equation}
which satisfies (Theorem \ref{bdy0}) 
$ \partial_0 ({\bf w}_0) = {{\bf z}_1}'$,  
we have :
\begin{equation*}
1 
=<{ch_{*}}({\bf w}_0 ), [1]> 
= [1] {\hat{\otimes}}_{A} {\bf w}_0
= \delta_1  [V] {\hat{\otimes}}_{A} {\bf w}_0
= ( [V] {\hat{\otimes}}_{A \times_{\alpha} {\bf Z}} {\bf x} ) {\hat{\otimes}}_{A} {\bf w}_0 
\end{equation*}
\begin{equation}
= [V] {\hat{\otimes}}_{A \times_{\alpha} {\bf Z}} ( {\bf x} {\hat{\otimes}}_{A} {\bf w}_0 ) 
= [V] {\hat{\otimes}}_{A \times_{\alpha} {\bf Z} } {\partial_0  ( {\bf w}_0 )} 
= < {ch_{*}}( {{\bf z}_1}' ), [V] >
= 1.
\end{equation}  
In the same way,
\begin{equation*}
0 
= < {ch_{*}}( {\bf w}_0 ), \delta_1  [U] > 
= \delta_1  [U] {\hat{\otimes}}_{A} {\bf w}_0 
= [U] {\hat{\otimes}}_{A \times_{\alpha} {\bf Z} } {\partial_0 ( {\bf w}_0 )}
\end{equation*} 
\begin{equation}
= [U] {\hat{\otimes}}_{A \times_{\alpha} {\bf Z} } {{\bf z}_1}' 
= < {ch_{*}}({{\bf z}_1}' ), [U]> 
=0.
\end{equation}
We also saw that the even K-homology of 
${A_{\theta}}$ is   
${KK^0}({A_{\theta}}, {\bf C}) \cong {\bf Z}^2 $.
One generator,
${\bf Dirac}$, 
was described previously (\ref{Athetadirac}). 

\begin{Prop}
\label{partialonewoneisdirac}
Under the map 
$\partial_1 : {KK^1}(A,{\bf C}) \rightarrow {KK^0}({A_{\theta}},{\bf C})$
we have\\
$ \partial_1  ({\bf w}_1) = {\bf Dirac}$, as elements of K-homology. 
\end{Prop}
\begin{proof}
 We recall that $ \partial_1  ({\bf w}_1) = {\bf x} {\hat{\otimes}}_{A \hat{\otimes} {\bf C}_1} {\bf w}_1$,
 where ${\bf x} \in {KK^1}( A_{\theta}, A)$ is the Pimsner-Voiculescu 
 Toeplitz element, described in (\ref{defntoeplitzelement}). 
The Toeplitz element ${\bf x}$  is represented by the Kasparov triple
\begin{equation}
\label{toeplitzagain}
(E_1,\phi_1,F_1) =
 ({l^2}({\bf Z},A) \hat{\otimes} {\bf C}_1,
{\phi_1}' \hat{\otimes} 1,
  F \hat{\otimes} \epsilon)
 \in
{\bf E}( A_{\theta} , A \hat{\otimes} {\bf C}_1 )
\end{equation}
 with 
\begin{equation}
\phi_1 (a) ( \xi \hat{\otimes} \omega) =
{\phi_1}' (a) \xi \hat{\otimes} \omega 
\end{equation}
 where ${\phi_1}' : A_{\theta}$
 $\rightarrow {\bf B} ( {l^2}({\bf Z},A)  )$ is defined by
\begin{equation}
({\phi_1}' (U) \xi)(n) = {\lambda^{-n}} U \xi(n), \quad
({\phi_1}' (V) \xi)(n) = \xi(n+1),
\end{equation}
 (remember that $A =$ ${C^{*}}(U)$)
 and $(F \xi) (n) = sign (n) \xi(n)$.

The Fredholm module ${\bf w}_1$ 
$= ( H = {l^2}({\bf Z}), \pi, F)$ 
 is represented by the Kasparov triple
\begin{equation}
{\bf w}_1 =
\left(
 H  \hat{\otimes} {\bf C}_1,  \pi \hat{\otimes} 1,
 F \hat{\otimes} \epsilon
\right)
 \in
{\bf E}( A, {\bf C}_1 )
\end{equation}
 where 
$(\pi \hat{\otimes} 1)(a) (\xi \hat{\otimes} \omega) =$ 
$\pi(a) \xi \hat{\otimes} \omega$, 
 and $(F \xi)(n) = sign(n) \xi(n)$, for $\xi \in H$, $\omega \in {\bf C}_1$. 
 To take the Kasparov product we need the Kasparov triple
 $\tau_{{\bf C}_1} ({\bf w}_1 ) \in$
 ${\bf E}( A \hat{\otimes} {\bf C}_1, {\bf C}_1 \hat{\otimes} {\bf C}_1)$.
 We have from (\ref{defntaumap}) that:
\begin{equation}
{\tau_{{\bf C}_1}} ({\bf w}_1 ) =
\left(
H \hat{\otimes} ({\bf C}_1 \hat{\otimes} {\bf C}_1),
 \pi \hat{\otimes} ( 1 \hat{\otimes} id),
 F  \hat{\otimes} ( \epsilon \hat{\otimes} 1)
\right)
 \in {\bf E}( A \hat{\otimes} {\bf C}_1 , {\bf C}_1 \hat{\otimes} {\bf C}_1).
\end{equation}
 We can identify ${\bf C}_1 \hat{\otimes} {\bf C}_1$
 with $M_2 ( {\bf C})$, via the map
\begin{equation}
\epsilon \hat{\otimes} 1 \mapsto
\left(
\begin{array}{cc}
0 & 1 \cr
1 & 0 \cr
\end{array}
\right),
\quad
1 \hat{\otimes} \epsilon \mapsto
\left(
\begin{array}{cc}
0 & -i \cr
i & 0 \cr
\end{array}
\right)
\end{equation}
 and realise
$\tau_{{\bf C}_1} ({\bf w}_1 )$
 as an element of
${\bf E}( A \hat{\otimes} {\bf C}_1 , M_2({\bf C}))$.
 This identification gives 
\begin{equation}
\label{defntauconewone}
{\tau_{{\bf C}_1}} ({\bf w}_1 ) 
\cong
(
 H \hat{\otimes} M_2 ({\bf C}),
 \pi \hat{\otimes} \rho,
  F \hat{\otimes}
\left(
\begin{array}{cc}
0 & 1 \cr
1 & 0 \cr
\end{array}
\right)
)
 =
 (E_2, \phi_2, F_2) 
 \in {\bf E}( A \hat{\otimes} {\bf C}_1 , M_2({\bf C}))
\end{equation}
 where
\begin{equation}
(\pi \hat{\otimes} \rho)(a \hat{\otimes} 1) ( \xi \hat{\otimes} T) =
\pi(a) \xi \hat{\otimes} T,
\end{equation}
\begin{equation}
(\pi \hat{\otimes} \rho)(a \hat{\otimes} \epsilon) ( \xi \hat{\otimes} T) =
\pi(a) \xi \hat{\otimes}
 \left(
\begin{array}{cc}
0 & -i \cr
i & 0 \cr
\end{array}
\right)
T.
\end{equation}
 The Kasparov product 
${\bf x} {\hat{\otimes}}_{A \hat{\otimes} {\bf C}_1} \tau_{ {\bf C}_1} ({\bf w}_1 )$
of the two triples (\ref{toeplitzagain}), (\ref{defntauconewone})
 is given  \cite{blackadar}, p166, by  the triple 
\begin{equation}
\label{firstproduct}
(E = E_1 {\hat{\otimes}}_{\phi_2} E_2 ,
 \phi = \phi_1 \hat{\otimes} 1, F)  
 \in {\bf E}( A_{\theta}, M_2( {\bf C}))
\end{equation}
 where the difficulty lies in finding a suitable $F$. 
 We have 
\begin{equation}
E_1 {\hat{\otimes}}_{\phi_2} E_2 \cong
 ( {l^2}({\bf Z},A) \hat{\otimes} {\bf C}_1 )
 {\hat{\otimes}}_{\phi_2}
 ( H \hat{\otimes} M_2 ({\bf C}) ),
\end{equation}
with
\begin{equation}
( \delta_k a_k \hat{\otimes} 1) \hat{\otimes} ( \xi \hat{\otimes} T)
\sim
( \delta_k  \hat{\otimes} 1) \hat{\otimes} ( \pi(a_k) \xi \hat{\otimes} T),
\end{equation}
\begin{equation}
( \delta_k a_k \hat{\otimes} \epsilon ) \hat{\otimes} ( \xi \hat{\otimes} T)
\sim
( \delta_k  \hat{\otimes} 1 ) \hat{\otimes}
( \pi(a_k) \xi \hat{\otimes}
\left(
\begin{array}{cc}
0 & -i \cr
i & 0 \cr
\end{array}
\right)
T).
\end{equation}
 Hence we can identify
$( {l^2}({\bf Z},A) \hat{\otimes} {\bf C}_1 )
 {\hat{\otimes}}_{\phi_2}
 ( H \hat{\otimes} M_2 ({\bf C}) )$
 as a submodule of
 ${l^2}({\bf Z}^2) \hat{\otimes} M_2 ({\bf C})$
 via the map 
\begin{equation}
( \delta_k a_k \hat{\otimes} 1) \hat{\otimes} ( e_l \hat{\otimes} T)
 \mapsto
( \delta_k \otimes \pi(a_k) e_l)
\hat{\otimes} T.
\end{equation}
 Under these identifications, we have that
$\phi = \phi_1 \hat{\otimes} 1$ acts by
\begin{equation}
\phi(U)
(( \delta_k \otimes  e_l)
\hat{\otimes} T)
=
( {\lambda^{-k}} \delta_k \otimes  e_{l+1})
\hat{\otimes} T,
\end{equation}
\begin{equation}
\phi(V)
(( \delta_k \otimes  e_l)
\hat{\otimes} T)
=
( \delta_{k+1} \otimes  e_l)
\hat{\otimes} T.
\end{equation}

We calculate the operator $F$ for the product via the Connes-Skandalis formalism of connections \cite{blackadar}, p170. 
 We know that there exists an 
$F_2$-connection $G$ for $E_1$, and from (\ref{suitableF}), having found such a $G$, an appropriate $F$ for the product is given by 
\begin{equation}
F = F_1 \hat{\otimes} 1 + ( (1 - {F_1}^2 )^{1/2} \hat{\otimes} 1) G.
\end{equation} 
 Since ${F_1}^2 =1$ in our situation, we can take 
$F = F_1 \hat{\otimes} 1$. There is no need to explicitly find $G$, knowledge of its existence is enough. 

We also have
 $F_1 \hat{\otimes} 1$ acting on $E$ (as a submodule of 
 ${l^2}({\bf Z}^2) \hat{\otimes} M_2({\bf C})$ ) by 
\begin{equation}
(F_1 \hat{\otimes} 1)( ( \delta_k \otimes  e_l)
\hat{\otimes} T)
=
(sign(k) \delta_k \otimes  e_l)
\hat{\otimes}
\left(
\begin{array}{cc}
0 & -i \cr
i & 0 \cr
\end{array}
\right)
T.
\end{equation}
 Note that the submodule that we have identified with $E$ is invariant under the action of $A_{\theta}$ and of $F$. 

 We have calculated the product as a Kasparov triple in 
${\bf E}( A_{\theta} , M_2 ({\bf C}))$. 
 We need to show that this represents the same element of K-homology as the Fredholm module ${\bf Dirac}$. 
 We will exhibit a homotopy of Kasparov triples from (\ref{firstproduct}) to a new element ${\bf y}_1$ of 
 ${\bf E}( A_{\theta} , M_2 ({\bf C}))$. 
 Then we use the KK-equivalence of $M_2 ({\bf C})$ and ${\bf C}$ to obtain a Kasparov triple in ${\bf E}( A_{\theta} ,{\bf C})$, still representing the product, which also represents ${\bf Dirac}$. 

 The homotopy of Kasparov triples ${ \{ {\bf y}_t \} }_{0 \leq t \leq 1}$ is given by : 
\begin{equation}
{\bf y}_t =
( E, \phi, {F_t}' ) 
\in {\bf E}( A_{\theta} , M_2 ({\bf C}))
\end{equation}
 with 
\begin{equation}
{F_t}'  
( (\delta_k \otimes e_l ) \hat{\otimes} T) =
 ( k^2 + t^2 l^2 )^{-1/2} 
(\delta_k \otimes e_l ) \hat{\otimes} 
\left(
\begin{array}{cc}
0 & ik + tl \cr
-ik + tl & 0 \cr
\end{array}
\right)
T.
\end{equation}
Now, ${\bf y}_0$ is the triple representing the product 
 $\partial_1 ({\bf w}_1)$ 
calculated in (\ref{firstproduct}), while ${\bf y}_1$ is the triple 
\begin{equation}
 (E, \phi, {F_1}')  \in 
{\bf E}( A_{\theta}, M_2 ( {\bf C}))
\end{equation}
 which therefore also represents the product.

We now use the KK-equivalence of
$M_2 ({\bf C})$ and ${\bf C}$ to realise the product as an element of 
${\bf E}( A_{\theta},{\bf C})$.
 The KK-equivalence is implemented (on the right) by the Kasparov triple
\begin{equation}
{\bf z}= ( {\bf C}^2, id,0) \in {\bf E}( M_2 ({\bf C}), {\bf C}).
\end{equation}
 Taking the product with ${\bf z}$ gives us the triple
\begin{equation}
(E' = E {\hat{\otimes}}_{id} {\bf C}^2 , \phi' = \phi \hat{\otimes} 1,
F' = {F_1}' \hat{\otimes} 1 )
 \in {\bf E}( A_{\theta}, {\bf C}).
\end{equation}
 The same argument as above tells us that this is the appropriate $F'$ for the product. 
 We can identify $E' =$ 
 $( {l^2}({\bf Z}^2) \hat{\otimes} M_2({\bf C})) {\hat{\otimes}}_{id} {\bf C}^2$
 with
${l^2}({\bf Z}^2) \oplus {l^2}({\bf Z}^2)$
  via the map 
\begin{equation}
(
\xi
\hat{\otimes} I_2)
\hat{\otimes}
\left[
\begin{array}{c}
\alpha \cr
\beta \cr
\end{array}
\right]
\mapsto
\left[
\begin{array}{c}
\alpha \xi \cr
\beta \xi \cr
\end{array}
\right].
\end{equation}
Under this identification, we have $A_{\theta}$ acting via
$\phi' = \pi' \oplus \pi'$, with
\begin{equation}
{\pi'}(U) ( \delta_k \otimes e_l) =   {\lambda^{-k}} \delta_k \otimes e_{l+1}, \quad
{\pi'}(V) ( \delta_k \otimes e_l) =  \delta_{k+1} \otimes e_l,
\end{equation}
 and 
$F' = 
\left(
\begin{array}{cc}
0 & F'' \cr
{F''}^{*} & 0 \cr
\end{array}
\right)$, 
 where 
\begin{equation}
F'' ( \delta_k \otimes e_l) =
\left\{
\begin{array}{cc}
{\frac{ik+l}{(k^2+l^2)^{1/2}}} \delta_k \otimes e_l & : (k,l) \neq (0,0) \cr
\delta_0 \otimes e_0  & : (k,l) = (0,0) \cr
\end{array}
\right.
\end{equation}

Now, let $\{ e_{p,q} \}$ be an arbitrary  new orthonormal basis for ${l^2}({\bf Z}^2)$, and define a unitary operator
$Q : {l^2}({\bf Z}^2) \rightarrow {l^2}({\bf Z}^2)$ by
\begin{equation}
Q ( \delta_k \otimes e_l) = {\lambda}^{kl} e_{l,k}.
\end{equation}
 Then the triple $(E', \phi', F')$ is unitarily equivalent to the triple 
 $(E, \phi, F) \in$ 
 ${\bf E} (A_{\theta}, {\bf C})$, where
\begin{equation}
 E = E' = {l^2}({\bf Z}^2) \oplus {l^2}({\bf Z}^2),
\end{equation}
 with $\phi = \pi \oplus \pi$ acting via
\begin{equation}
\pi(U) e_{p,q} = e_{p+1,q}, \quad
\pi(V) e_{p,q} = {\lambda^{p}} e_{p,q+1},
\end{equation}
and 
\begin{equation}
F =
\left(
\begin{array}{cc}
0 & {F_0}^{*} \cr
{F_0} & 0 \cr
\end{array}
\right), \quad
{F_0} e_{p,q}
=
\left\{
\begin{array}{cc}
{\frac{p+iq}{(p^2 +q^2)^{1/2}}} e_{p,q} & : (p,q) \neq (0,0) \cr
e_{0,0} & : (p,q) = (0,0) \cr
\end{array}
\right.
\end{equation}
We recognise this as a triple representing the Fredholm module ${\bf Dirac}$,
as defined in (\ref{Athetadirac}).
 So we have shown by direct computation  that
${\partial_1}({\bf w}_1) = {\bf Dirac}$ as elements of K-homology.
\end{proof}

\begin{Lemma} 
\label{istardiraciszero}
${i^{*}}({\bf Dirac}) = {\bf 0} \in {KK^0}(A,{\bf C})$. 
\end{Lemma}
\begin{proof}
 We have  ${i^{*}}({\bf Dirac}) \in {KK^0}(A,{\bf C})$ is the Fredholm module 
\begin{equation}
\label{defnistardirac}
(
{l^2}( {\bf Z}^2) \oplus {l^2}( {\bf Z}^2), {\pi_0} \oplus {\pi_0}, 
F = \left(
\begin{array}{cc}
0 & {F_0}\cr
{F_0} & 0 \cr
\end{array}
\right)
) 
\end{equation}
 where ${\pi_0}(U) e_{p,q} = e_{p+1,q}$, 
and 
\begin{equation}
{F_0}{e_{p,q}} = \left\{
\begin{array}{cc}
{\frac{p+iq}{(p^2 + q^2)^{1/2}}} e_{p,q} & : (p,q) \neq (0,0) \cr
e_{0,0} & : (p,q) =(0,0) \cr
\end{array}
\right.
 \end{equation}
 Now, ${K_0}(A) \cong {\bf Z}$, generated by $[1]$, 
 ${KK^0}(A,{\bf C}) \cong {\bf Z}$, generated by ${\bf w}_0$, 
 and
 $<{ch_{*}}({\bf w}_0),[1]>$ $ =1$. 
 If ${i^{*}}({\bf Dirac})$ is a nontrivial element of K-homology, then we will have ${i^{*}}({\bf Dirac}) = n{\bf w}_0$, for some $n \neq 0$, and so 
 $< {ch_{*}}({i^{*}}({\bf Dirac})), [1]> =n$.
 But we see from (\ref{defnistardirac}) 
that $[F, ({\pi_0} \oplus {\pi_0})(1) ] =0$, hence 
  $< {ch_{*}}({i^{*}}({\bf Dirac})), [1]> =0$, and so 
 ${i^{*}}({\bf Dirac})$ represents a trivial element of K-homology. 
\end{proof}

We calculate, using ${\partial_1}({\bf w}_1 ) = {\bf Dirac}$, that
\begin{equation*}
1
= <{ch_{*}}({\bf w}_1),[U]>
=[U]{\hat{\otimes}}_{A} {\bf w}_1
=\delta_0  [p] {\hat{\otimes}}_{A} {\bf w}_1 
=[p] {\hat{\otimes}}_{A \times_{\alpha} {\bf Z}} {\bf Dirac}
= <{ch_{*}}({\bf Dirac}),[p]>.
\end{equation*}
Also,
\begin{equation*}
0 
= <{ch_{*}}( {\bf Dirac} ),[1]>
= [1]{\hat{\otimes}}_{A \times_{\alpha} {\bf Z}} {\bf Dirac}
= [1]{\hat{\otimes}}_{A \times_{\alpha} {\bf Z}} \partial_1  ({\bf w}_1) 
= \delta_0  [1] {\hat{\otimes}}_{A} {\bf w}_1 
=0.
\end{equation*}
 We want to describe a second generator 
of
${KK^0}(A_{\theta},{\bf C}) \cong {\bf Z}^2$.
 We will denote this generator by ${\bf z}_0$. 

It is easy to see that the map 
$(id - {\alpha^{*}}) : {KK^0}(A,{\bf C}) \rightarrow {KK^0}(A,{\bf C})$
 is the zero map.
Hence  the map $i^{*}$ is surjective, and 
 ${i^{*}}({\bf Dirac}) ={\bf 0}$, we impose for all values of $\theta$ that  
 ${i^{*}}({\bf z}_0) = {\bf w}_0 \in {KK^0}(A,{\bf C})$.
 However, it is difficult to describe such a Fredholm module explicitly. 
 In the case 
$\theta=0$,  we take ${\bf z}_0$ 
 to be the canonical even Fredholm module 
 (Lemma \ref{evenfred}) 
over 
$C({\bf T}^2)$, 
 and for other values of $\theta$ we want the corresponding ${\bf z}_0$ to be a ``continuous deformation'' of this.

Under the map
${i_{*}} : {K_0}(A) \rightarrow {K_0}(A_{\theta})$,
we have
${i_{*}}[1]=[1]$,
and we calculate that
\begin{equation*}
1
=<{ch_{*}}({\bf w}_0 ),[1]>
= [1] {\hat{\otimes}}_{A} {\bf w}_0
= [1] {\hat{\otimes}}_{A} {i^{*}}({\bf z}_0) 
= {i_{*}} [1] {\hat{\otimes}}_{A \times_{\alpha} {\bf Z} } {\bf z}_0
= <{ch_{*}}({\bf z}_0 ),[1]>.
\end{equation*}
We also need to know 
$<{ch_{*}}({\bf z}_0 ),[p]>$. 
Since in the case $\theta = 0$ the Powers-Rieffel projection $p$ \cite{davidson}, p170, is just $p=0$, and we want our ${\bf z}_0$ to be a deformation of the $\theta=0$ case, we will require that 
\begin{equation}
\label{zzerorieffelproj}
 <{ch_{*}}({\bf z}_0 ),[p]>=0. 
\end{equation}

 In the case $\theta$ is rational, 
$\theta = m/n$, with $m$, $n$ relatively prime integers, $n >0$, 
then $A_{\theta}$ is the algebra of continuous sections of a vector bundle over ${\bf T}^2$, whose fibres are full matrix algebras $M_n ({\bf C})$. 

We construct a Fredholm module ${{\bf z}_0}'$ over $A_{\theta}$ as follows:
\begin{equation}
\label{zzeroprime}
{{\bf z}_0}' = (
{\bf C}^n \oplus {\bf C}^n, 
 {\pi_0}' = \varphi \oplus 0, 
 {F_0}' = 
\left(
\begin{array}{cc}
0 & 1 \cr
1 & 0 \cr
\end{array}
\right)
)
\end{equation}
 with 
\begin{equation}
\varphi(U) =
\left(
\begin{array}{cccc}
0 & & & 1 \cr
1 &  & & .. \cr
 & & 0 & 0 \cr
0 & .. & 1 & 0 \cr
\end{array}
\right), \quad
\varphi(V) =
\left(
\begin{array}{cccc}
1 & & &  \cr
  & \lambda  & &  \cr
 & & .. &  \cr
 & &  & \lambda^{n-1} \cr
\end{array}
\right).
\end{equation}
 Then an easy calculation (see Lemma \ref{canonicalevenpairing}) shows that 
\begin{equation}
< ch_{*} ( {{\bf z}_0}' ), [1]> = n,
\end{equation}
 and further, under ${i^{*}}: {KK^0}(A_{\theta}, {\bf C}) \rightarrow 
{KK^0}(A,{\bf C})$, 
 we have ${i^{*}}( {{\bf z}_0}') = n {\bf w}_0$. 
 Since $i^{*}$ is surjective, ${{\bf z}_0}'$ cannot be a generator of 
 ${KK^0}(A_{\theta}, {\bf C})$. In particular, 
 ${{\bf z}_0}' \neq {\bf z}_0$.

It would be good to have an explicit description of the Fredholm module ${\bf z}_0$. 
 In the final  section of this paper, we describe an approach to this via 
Pimsner and Voiculescu's work on embedding the irrational rotation algebras in AF-algebras \cite{pv80a}. 
 The general question of concrete realizations (as Fredholm modules) of the K-homology of AF-algebras is not well-studied and is an interesting topic for future research.

\section{K-homology of the $A_{\theta}$ via AF-embeddings}

We conclude our study of the K-homology of the rotation algebras by  exploiting the AF embedding technique of Pimsner and Voiculescu \cite{pv80a} to try to find the missing generator of the even K-homology. 

Given an irrational $\theta \in (0,1)$, Pimsner and Voiculescu constructed an embedding of the irrational rotation algebra $A_{\theta}$ in an AF-algebra $C_{\theta}$ as follows. 
We begin by  considering  the continued fraction expansion
\begin{equation}
\theta = {lim_{n \rightarrow \infty}} [ a_0 , a_1 , ..., a_n ]
= {lim_{n \rightarrow \infty}} 
(a_0 + {\frac{1}{a_1 + {\frac{1}{ .. + {\frac{1}{a_n}}}}}})
\end{equation}
where $a_0 \in {\bf Z}$, and $a_1, .. a_n \in {\bf N}$. 
The rational approximations 
${\frac{p_n}{q_n}} = [a_0 , ... , a_n]$ are given recursively by 
\begin{equation*}
p_0 = a_0 , \,\, q_0 = 1, \,\, p_1 = a_0 a_1 +1, \,\, q_1 = a_1,
\end{equation*}
\begin{equation}
\label{pnqn}
p_n = a_n p_{n-1} + p_{n-2}, \,\, q_n = a_n q_{n-1} + q_{n-2}.
\end{equation}
 We define a sequence ${\{ C_n \}}_{n \geq 1}$ of finite-dimensional C*-algebras by 
\begin{equation}
C_n = M_{q_n} ( {\bf C}) \oplus M_{q_{n-1}}({\bf C}). 
\end{equation}
The maps
$\phi_{n+1,n} : C_n \rightarrow C_{n+1}$ are given by 
\begin{equation}
\label{phinnplusone}
\phi_{n+1,n} : 
\left(
\begin{array}{cc}
A_n & 0 \cr
0 & B_n \cr
\end{array}
\right)
\mapsto 
\left(
\begin{array}{cc}
{W_{n+1}} & 0 \cr
0 & I \cr
\end{array}
\right)
\left(
\begin{array}{ccc}
\begin{array}{ccc}
A_n & .. & 0 \cr
.. & .. & 0 \cr
0 & ..& A_n \cr
\end{array}
& 0 & 0\cr
0 & B_n & 0 \cr
0 & 0 & A_n \cr
\end{array}
\right)
\left(
\begin{array}{cc}
W_{n+1}^{*} & 0 \cr
0 & I \cr
\end{array}
\right)
\end{equation}
 where the $A_n$'s in the top left corner occur with multiplicity $a_n$, and $W_{n+1} \in {M_{q_{n+1}}}( {\bf C}) $ is a unitary.

The AF-algebra $C_{\theta}$ is defined to be 
$C_{\theta} = {\lim_{n \rightarrow \infty}} C_n$. 
 For each $n$ we have an inclusion map $\phi_n : C_n \rightarrow C_{\theta}$.
 Pimsner and Voiculescu \cite{pv80a} proved the following result :

\begin{Thm}
 There is an injective *-homomorphism $\rho : A_{\theta} \rightarrow C_{\theta}$, such that 
${\rho_{*}} : {K_0}( A_{\theta} ) \rightarrow {K_0}(C_{\theta}) $ is an isomorphism of abelian groups. 
 Furthermore, if $\tau$ and $\sigma$ are the canonical normalized traces on $A_{\theta}$ and $C_{\theta}$ respectively, then 
$ \tau_{*} = \sigma_{*} \rho_{*}$ is an order isomorphism of ${K_0}(A_{\theta})$ onto ${\bf Z} + {\bf Z} \theta$.  
\end{Thm}

So ${K_0}(C_{\theta}) \cong {\bf Z}^2$, generated by 
$[\rho(1)]$ and $[\rho(p)]$, where $[1]$ and $[p]$ generate ${K_0}(A_{\theta})$, 
and further ${K_1}(C_{\theta}) \cong 0$ (since $C_{\theta}$ is AF). 
 Hence by the universal coefficient theorem we have 
${KK^0}(C_{\theta}, {\bf C}) \cong {\bf Z}^2$. 
(We will also calculate this directly, in a way that will be more useful for our purposes.)
 Since 
$C_{\theta}$ 
is AF, for each 
$[x] \in {K_0}(C_{\theta})$ 
there exists $n$, and $[x_n] \in {K_0}(C_n)$, such that 
$[ \phi_n (x_n) ] = [x]$. 
 In fact, under 
$\phi_1 : C_1 \rightarrow C_{\theta}$, 
we have :
\begin{Lemma}
$[ \phi_1 (1)] = [ \rho(1)]$, and
$[\phi_1 (p_1) ] = [ \rho(p) ]$, 
 where 1 is the unit of $C_1$, and 
 $p_1$ is the rank one projection in $C_1 \cong$ ${M_{a_1}}({\bf C}) \oplus {\bf C}$ given by 
\begin{equation}
p_1 = 
\left(
\begin{array}{cccc}
1 & 0 & ..& 0 \cr
0 & 0 & ..& 0 \cr
.. & .. & .. & .. \cr
0 & ..& .. & 0 \cr
\end{array}
\right)
\end{equation}
 \end{Lemma}

 We have 
\begin{equation}
C_1 \rightarrow C_2 \rightarrow ... \rightarrow C_n \rightarrow ... \rightarrow C_{\theta} \leftarrow^{\rho} A_{\theta} 
\end{equation}
 If we could exhibit the Fredholm modules generating the even K-homology of $C_{\theta}$, then we could pull these back via $\rho$ to get Fredholm modules over $A_{\theta}$. 

 For each $n$, we have ${KK^0}(C_n, {\bf C}) \cong {\bf Z}^2$, generated by Fredholm modules ${\bf z}_1^{(n)}$ and ${\bf z}_2^{(n)}$. 
These are defined as follows (see Example \ref{evenfred}) :
\begin{equation}
{\bf z}_1^{(n)} = 
( {\bf C}^{q_n} \oplus {\bf C}^{q_n}, 
\left(
\begin{array}{cc}
{A_n} & 0 \cr
0 & {B_n} \cr
\end{array}
\right)
 \mapsto 
\left(
\begin{array}{cc}
{A_n} & 0 \cr
0 & 0 \cr
\end{array}
\right), 
F = 
\left(
\begin{array}{cc}
0 & 1 \cr
1 & 0 \cr
\end{array}
\right)
)
\end{equation}
\begin{equation}
{\bf z}_2^{(n)} =
( {\bf C}^{q_{n-1}} \oplus {\bf C}^{q_{n-1}}, 
\left(
\begin{array}{cc}
{A_n} & 0 \cr
0 & {B_n} \cr
\end{array}
\right)
 \mapsto 
\left(
\begin{array}{cc}
{B_n} & 0 \cr
0 & 0 \cr
\end{array}
\right), 
F = 
\left(
\begin{array}{cc}
0 & 1 \cr
1 & 0 \cr
\end{array}
\right)
)
\end{equation}
 Under the maps $\phi_{n+1,n} : C_n \rightarrow C_{n+1}$ (\ref{phinnplusone})
 we have 
\begin{equation}
{\phi_{n+1,n}^{*}} ( {\bf z}_1^{(n+1)}) = a_n {\bf z}_1^{(n)}  + 
{\bf z}_2^{(n)},
\end{equation}
\begin{equation}
{\phi_{n+1,n}^{*}} ( {\bf z}_2^{(n+1)}) =   {\bf z}_1^{(n)}.
\end{equation}
 Hence all the maps 
\begin{equation}
{\phi_{n+1,n}^{*}} : {KK^0}( C_{n+1}, {\bf C}) \cong {\bf Z}^2 \rightarrow 
{KK^0}( C_{n}, {\bf C}) \cong {\bf Z}^2,
\end{equation}
 are surjective, because the matrix 
$\left(
\begin{array}{cc}
a_n & 1 \cr
1 & 0 \cr
\end{array}
\right)$
is always invertible.  

 We now appeal to the following special case of a much more general result of Rosenberg and Schochet \cite{rs}.

\begin{Prop}
Suppose that $A = {lim_{\rightarrow}} A_n$ is an AF-algebra. Then the following sequences on K-homology are exact:
\begin{equation*}
0 \rightarrow {lim^1}_{\leftarrow}  {KK^1}( A_n, {\bf C}) \rightarrow {KK^0}(A,{\bf C}) \rightarrow {lim_{\leftarrow}} {KK^0}( A_n , {\bf C}) \rightarrow 0
\end{equation*}
\begin{equation*}
0 \rightarrow {lim^1}_{\leftarrow}  {KK^0}( A_n, {\bf C}) \rightarrow {KK^1}(A,{\bf C}) \rightarrow {lim_{\leftarrow}} {KK^1}( A_n , {\bf C}) \rightarrow 0
\end{equation*}
\end{Prop}

The left hand term is Milnor's ${lim^1}_{\leftarrow}$, and the right hand term is the inverse limit of the K-homology groups. 

 It follows from \cite{weibel}, p80 that if the maps 
${KK^i}(A_{n+1}, {\bf C}) \rightarrow {KK^i}(A_n, {\bf C})$ are all surjective, then both $lim^1$ terms vanish, and hence 
${KK^i}(A,{\bf C}) \cong$ 
 ${lim_{\leftarrow}} {KK^i}( A_n , {\bf C})$, $(i=0,1)$. 
 This is true in our situation, and furthermore we have 
${lim_{\leftarrow}} {KK^0}( C_n , {\bf C})$  $\cong {\bf Z}^2$, 
${lim_{\leftarrow}} {KK^1}( C_n , {\bf C})$  $\cong 0$.
 Hence ${KK^0}(C_{\theta},{\bf C})$  $\cong {\bf Z}^2$, and 
${KK^1}(C_{\theta},{\bf C})$  $\cong 0$, which we knew already, but via the definition of the inverse limit we can now visualize the elements of ${KK^0}( C_{\theta}, {\bf C})$. 

\begin{Lemma}
An element ${\bf z} \in {KK^0}( C_{\theta}, {\bf C})$ is represented by a sequence of Fredholm modules ${\{ {\bf z}_n \}}_{n \geq 1}$, with 
${\bf z}_n \in {KK^0}(C_n, {\bf C})$, such that under each of the inclusion maps 
$\phi_n : C_n \rightarrow C_{\theta}$, we have 
${\phi_n^{*}}({\bf z}) = {\bf z}_n$. 
 It follows immediately that for each of the maps 
 $\phi_{n+k,n} : C_n \rightarrow C_{n+k}$, we have 
${\phi_{n+k,n}^{*}} ({\bf z}_{n+k}) = {\bf z}_n$. 
\end{Lemma}

\begin{proof} 
This follows immediately from the definition of the inverse limit of a sequence of abelian groups. 
Recall that, for 
a tower of abelian groups 
\begin{equation}
 \ldots \rightarrow G_n \rightarrow^{f_n} G_{n-1} \rightarrow^{f_{n-1}} \ldots \rightarrow^{f_0} G_0
\end{equation}
 the inverse limit ${\lim_{\leftarrow}}{G_n}$ is isomorphic to the abelian group consisting of all sequences ${\{ g_n \}}_{n \geq 0}$, with $g_n \in G_n$ for each $n$, such that $g_{n-1} = {f_n}(g_n)$. 
\end{proof}

Furthermore, given any $[x] \in {K_0}(C_{\theta})$, there exists $n$, and $[x_n] \in {K_0}(C_n)$, such that 
$[ \phi_n (x_n) ] = [x]$. 
 So 
\begin{equation}
< {ch_{*}}({\bf z}) , [x]> =
< {ch_{*}}({\bf z}) , [ \phi_n (x_n) ]> =
< {ch_{*}}( {\phi_n^{*}}({\bf z})) , [x_n]> =
< {ch_{*}}({\bf z}_n ) , [x_n]>.
\end{equation}

 We want to find a Fredholm module ${\bf z}_0 \in $ 
${KK^0}( C_{\theta}, {\bf C})$ 
so that 
$\rho^{*} ({\bf z}_0 ) = {\bf w}_0$ 
$ \in {KK^0}( A_{\theta}, {\bf C})$. 
 We need 
\begin{equation}
< {ch_{*}}({\bf z}_0) , [\rho(1)]> =
< {ch_{*}}({\bf w}_0), [1]> =1,
\end{equation}
 and 
\begin{equation}
< {ch_{*}}({\bf z}_0) , [\rho(p)]> =
< {ch_{*}}({\bf w}_0), [p]> =0,
\end{equation}
 We can  take
${\bf z}_1 = {\phi_1^{*}}( {\bf z}_0 ) =$
$ {\bf z}_2^{(1)}$, 
since 
$< {ch_{*}}({\bf z}_1^{(1)}), [1]> = a_1$, 
 $< {ch_{*}}({\bf z}_1^{(1)}), [p_1]>=1$,\\ 
 $< {ch_{*}}({\bf z}_2^{(1)}), [1]> = 1$ and 
 $< {ch_{*}}({\bf z}_2^{(1)}), [p_1]>=0$.

We can calculate the corresponding 
${\bf z}_n = {\phi_n^{*}}( {\bf z}_0 ) \in $ 
${KK^0}(C_n, {\bf C})$ 
in the same way, provided we know all the $a_n$'s. 
 We have 
\begin{equation}
{\bf z}_n = x {\bf z}_1^{(n)} + y {\bf z}_2^{(n)}
\end{equation}
where 
\begin{equation*}
\left(
\begin{array}{c}
x \cr
y \cr
\end{array}
\right)
=
{
\left(
\begin{array}{cc}
a_{n-1} & 1 \cr
1 & 0 \cr
\end{array}
\right)
}^{-1}
 \ldots
{
\left(
\begin{array}{cc}
a_1 & 1 \cr
1 & 0 \cr
\end{array}
\right)
}^{-1}
\left(
\begin{array}{c}
0 \cr
1 \cr
\end{array}
\right)
=
{
\left(
\begin{array}{cc}
q_n & q_{n-1} \cr
p_n & p_{n-1} \cr
\end{array}
\right)}^{-1}
\left(
\begin{array}{c}
0 \cr
1 \cr
\end{array}
\right)
\end{equation*}
\begin{equation} 
=
(-1)^{n}
\left(
\begin{array}{cc}
p_{n-1} & -q_{n-1} \cr
-p_n & q_n \cr
\end{array}
\right)
\left(
\begin{array}{c}
0 \cr
1 \cr
\end{array}
\right) 
 =
(-1)^{n}
\left(
\begin{array}{c}
-q_{n-1} \cr
q_n \cr
\end{array}
\right) 
\end{equation}
 (this follows from the relations (\ref{pnqn})). 
 
 In this way we obtain an element of 
 ${lim_{n \leftarrow \infty}}$  ${KK^0}( C_n, {\bf C})$ 
 representing ${\bf z}_0 \in$ 
 ${KK^0}(C_{\theta}, {\bf C})$. 
 It is not clear how to pull this back via $\rho$ to 
 ${KK^0}(A_{\theta}, {\bf C})$.

\section{Acknowledgements}
 I would like to thank my advisor, Professor Marc Rieffel, 
 for his advice and support throughout my time at Berkeley. 
 I am very grateful for his help. 
 I would also like to thank Erik Guentner, Nate Brown and Frederic 
 Latremoliere for many useful discussions.


\end{document}